\documentclass[12pt,reqno]{amsart}                  
\usepackage{amssymb}
\usepackage{mathrsfs}
\newtheorem{theorem}{Theorem}
\newtheorem{lemma}{Lemma}

\theoremstyle{remark}

\begin{document}
\title[Kreps-Yan theorem for Banach ideal spaces]{Kreps-Yan theorem for Banach ideal spaces}

\author{Dmitry B. Rokhlin}

\address{D.B. Rokhlin,
Faculty of Mathematics, Mechanics and Computer Sciences,
              Southern Federal University, 
Mil'chakova str., 8a, 344090, Rostov-on-Don, Russia}          
\email{rokhlin@math.rsu.ru}                   
\newcommand{\ri}{{\rm ri\,}}
\newcommand{\gr}{{\rm gr\,}}
\newcommand{\cl}{{\rm cl}}
\newcommand{\epi}{{\rm epi\,}}
\newcommand{\cone}{{\rm cone\,}}
\newcommand{\intern}{{\rm int\,}}
\newcommand{\conv}{{\rm conv\,}}
\newcommand{\dom}{{\rm dom\,}}
\newcommand{\supp}{{\rm supp\,}}
\newcommand{\esssup}{{\rm ess\,sup}}

\begin{abstract} 
Let $C$ be a closed convex cone in a Banach ideal space $X$ on a measurable space with a $\sigma$-finite measure. We prove that conditions $C\cap X_+=\{0\}$ and $C\supset -X_+$ imply the existence of a strictly positive continuous functional on $X$, whose restriction to $C$ is non-positive. 
\end{abstract}
\maketitle

Let $(\Omega,\mathscr F)$ be a measurable space, which is complete with respect to a measure (that is, a countably-additive function) $\mu:\mathscr F\mapsto [0,\infty]$. Consider the vector space $L^0(\mu)=L^0(\Omega,\mathscr F,\mu)$ of classes of $\mu$-equivalent a.s. finite $\mathscr F$-measurable functions. This space is a vector lattice (Riesz space) with respect to the natural order structure, induced by the cone $L^0_+(\mu)$ of non-negative elements
[1], [2]. Let $X$ be a solid subspace (ideal) in $L^0(\mu)$, that is, $X$ is a linear subset of $L^0(\mu)$ and conditions $x\in X$, $|y|\le |x|$ imply that $y\in X$. Assume that there is a norm on $X$, satisfying the condition $\|x\|\le\|y\|$, if $|x|\le |y|$, $x,y\in X$ (monotone norm) and $X$ is complete with respect to this norm. As is known, in this case $(X,\|\cdot\|)$ is called a \emph{Banach ideal space} on $(\Omega,\mathscr F,\mu)$ [1], [3].

Let $X$ be a Banach ideal space with the non-negative cone $X_+=\{x\in X:x\ge 0\}$. An element $g$ of the topological dual space $X'$ is called \emph{strictly positive} if $\langle x,g\rangle:=g(x)>0$, $x\in X_+\backslash\{0\}$. Slightly modifying the terminology of [4], we say that $X$ has the \emph{Kreps-Yan property}, if for any closed convex cone $C\subset X$, satisfying the conditions 
\begin{equation}
C\cap X_+=\{0\},\ \ \ -X_+\subset C,
\end{equation}
there exists a strictly positive element $g\in X'$ such that $\langle x,g\rangle\le 0$, $x\in C$. 

This property can be considered also in a more general setting of a locally convex space with a cone. We refer to [4] for the results in this direction as well as for further references and comments, concerning the papers of Kreps [5] and Yan [6].

It is said that a topological space $(X,\tau)$ has the \emph{Lindel\"of property}, if every open cover of $X$ has a
countable subcover. If the weak topology $\sigma(X,X')$ of a Banach space $X$ has the Lindel\"of property  (for brevity, $X$ is weakly Lindel\"of) and the set of strictly positive functionals $g\in X'$ is non-empty, then $X$ has the Kreps-Yan property [7, Theorem 1.1]. A space $X$ is weakly Lindel\"of if any of the following conditions hold true: (a) $X$ is reflexive, (b) $X$ is separable, (c) $X$ is weakly compactly generated (see [8]). 

Also, it is known that the space $L^\infty(\Omega,\mathscr F,\mathsf P)$ of $\mathscr F$-measurable essentially bounded functions (where $\mathsf P$ is a probability measure) has the Kreps-Yan property. This was established in [7, Theorem 2.1] (and, independently, in [9]). Note that the space $l^\infty$ of bounded sequences is not weakly Lindel\"of (see [10, Example 4.1(i)]). As long as $l^\infty$ can be considered as a Banach ideal space on the probability space $(\mathbb N,2^\mathbb N,\mathsf Q)$, $\mathsf Q(n)=1/2^n$, the Kreps-Yan property is strictly weaker then the Lindel\"of property of the weak topology.

On the other hand, as Example 2.1 of [4] shows, the space $l^1(\mathbb R_+)=\{(f_t)_{t\in\mathbb R_+}:\sum_{t\in\mathbb R_+}|f_t|<\infty\}$ does not possess the Kreps-Yan property, although the set of strictly positive functionals on $l^1(\mathbb R_+)$ is non-empty. Note that $l^1(\mathbb R_+)$ is a Banach ideal space on $(\mathbb R_+,2^{\mathbb R_+},\nu)$, where $\nu$ is the counting measure on $\mathbb R_+$. This measure is decomposable in the sense of the definition of Exercise 1.12.131 in [11], but is not $\sigma$-finite.

The mentioned example of [4] underlines the accuracy of our main result, which is formulated as follows.

\begin{theorem}
A Banach ideal space $X$ on $(\Omega,\mathscr F,\mu)$ with a $\sigma$-finite measure $\mu$ has the Kreps-Yan property.
\end{theorem}

It is natural to call this result the Kreps-Yan theorem for Banach ideal spaces.

Roughly speaking, the strategy of the proof follows the paper [7], where the argumentation was adapted to $L^\infty$. At first, we prove the existence of a strictly positive element $f\in X'$, which is bounded from above on the set $\{x\in C:\|x^-\|\le 1\}$ (Lemma 1). Then we establish the existence of an element $g\in X'$, $g\ge f$, whose restriction on $C$ is non-positive (Lemma 3). Here and in what follows, by $x^+$, $x^-$ we denote the positive and negative parts of an element $x$ of a vector lattice.

First of all, note that it is sufficient to consider the case of probability measure $\mu$. Actually, by the definition of $\sigma$-finitness there exists a disjoint sequence $(A_n)_{n=1}^\infty$, $A_n\in\mathscr F$, forming the partition of $\Omega$ and such that $\mu(A_n)<\infty$. The probability measure
$$ \mathsf P(A)=\sum_{n=1}^\infty\frac{\mu(A\cap A_n)}{2^n\mu(A_n)},\qquad A\in\mathscr F$$
is equivalent to $\mu$, that is, it has the same collection of negligible sets. Clearly, any Banach ideal space $X$ on $(\Omega,\mathscr F,\mu)$ can be regarded as a Banach ideal space on $(\Omega,\mathscr F,\mathsf P)$.

Denote by $\mathsf E x=\int x\,d\mathsf P$ the expectation with respect to the probability measure $\mathsf P$. 
Any Banach ideal space is a subset of the complete metric space $(L^0(\mathsf P),d)$, where
$$ d(x,y)=\mathsf E\frac{|x-y|}{1+|x-y|}.$$
Convergence in the topology, induced by this metric, coincides with the convergence in probability and boundedness coincides with the boundedness in probability. Recall that the set $H\subset L^0(\mathsf P)$ is said to be bounded in probability if for any $\varepsilon>0$ there exists $M>0$ such that $\mathsf P(|x|\ge M)<\varepsilon$, $x\in H.$

\begin{lemma} \label{lem:1}
Let $X$ be a Banach ideal space on $(\Omega,\mathscr F,\mathsf P)$ and let $C\subset X$ be a convex cone, which is closed in the norm topology of $X$ and satisfies conditions (1). Then there exists an equivalent to $\mathsf P$ probability measure $\mathsf Q$ such that
$$ \sup_{x\in C_1}\mathsf E_\mathsf Q x <\infty, \ \ \ C _1=\{x\in C:\|x^-\|\le 1\} $$
and $f=d\mathsf Q/d\mathsf P\in L^\infty(\mathsf P)$.
\end{lemma}
\textsc{Proof.} Consider the convex set
$$G=\conv |C_1|=\left\{y= \sum_{i=1}^m\alpha_i |x_i|:x_i\in C_1,\ \alpha_i\ge 0,\ \sum_{i=1}^m\alpha_i=1,\ m\in\mathbb N\right\}$$
and put
$$G^1=\{y\in L^1_+(\mathsf P): y\le x\ \textrm{for some}\ x\in G\},$$
where $L^1(\mathsf P)$ is the Banach lattice of $\mathsf P$-integrable functions. 
It is sufficient to prove that for any $A\in\mathscr F$ with $\mathsf P(A)>0$ there exists some $c\in\mathbb R_+$ such that 
\begin{eqnarray}
c I_A\not\in\overline{G^1-L^1_+(\mathsf P)},
\end{eqnarray}
where the overline denotes closure in $L^1(\mathsf P)$. Indeed, by Yan's theorem ([6, Theorem 2]) this implies the existence of an element $f\in L^\infty(\mathsf P)$ with $\mathsf P(f>0)=1$ such that $\sup_{y\in G^1}\mathsf E(y f)<\infty$. Note that  $y_n=x\wedge n\in G^1$ for any $x\in G$, $n\in\mathbb N$. By the monotone convergence theorem it follows that
$$\sup_{x\in C_1}\mathsf E_\mathsf Q x\le \sup_{x\in G}\mathsf E_\mathsf Q x=\sup_{y\in G^1}\mathsf E(y f)<\infty,$$
where $d\mathsf Q/d\mathsf P=f$.

To prove (2) assume the converse. Suppose that $A\in\mathscr F$, $\mathsf P(A)>0$ and $c I_A\in\overline{G^1-L^1_+(\mathsf P)}$ for all $c>0$. Then for any $\varepsilon>0$ there exists $z\in G$ such that 
\begin{eqnarray}
\mathsf P(A\cap\{z<1/\varepsilon\})<\varepsilon.
\end{eqnarray}
Actually, put $c>1/\varepsilon$, assume that the sequence $y_n\in G^1-L^1_+(\mathsf P)$ converges to $c I_A$ with probability $1$ and $y_n\le z_n$, $z_n\in G$. Then
$$ \mathsf P\left(\bigcap_{n=1}^\infty \bigcup_{m=n}^\infty (A \cap \{z_m<1/\varepsilon\})\right)\le 
\mathsf P\left(\bigcap_{n=1}^\infty \bigcup_{m=n}^\infty (A \cap \{y_m<1/\varepsilon\})\right)=0.$$
Consequently,
$$ \limsup_{n\to\infty}\mathsf P(A \cap \{z_n<1/\varepsilon\})\le \lim_{n\to\infty}\mathsf P\left(\bigcup_{m=n}^\infty (A \cap \{z_m<1/\varepsilon\})\right)=0.$$

An element $z\in G$, satisfying (3), admits the representation
$$ z=\sum_{i=1}^m\alpha_i |x_i|=\sum_{i=1}^m\alpha_i x_i+2\sum_{i=1}^m\alpha_i x_i^-=x+y,\ \ \ \alpha_i\ge 0,\  \sum_{i=1}^m\alpha_i=1, $$
where $x\in C_1$, $y\in X_+$, $\|y\|\le 2$. Using the inequality
$$\mathsf P(B_1\cap B_2)=\mathsf P(B_1)+\mathsf P(B_2)-\mathsf P(B_1\cup B_2)\ge\mathsf P(B_1)-\mathsf P(\Omega\backslash B_2),$$
we get
\begin{eqnarray}
\mathsf P\left(A\bigcap\left\{x<\frac{1}{2\varepsilon}\right\}\right) -
\mathsf P\left(y\ge\frac{1}{2\varepsilon}\right) &\le & \mathsf P\left(A\bigcap\left\{x<\frac{1}{2\varepsilon}\right\}\bigcap\left\{y<\frac{1}{2\varepsilon}\right\}\right)\nonumber\\
&\le & \mathsf P(A\cap\{z<1/\varepsilon\})<\varepsilon.
\end{eqnarray}

As is is known, the embedding of the Banach space $X$ in the metric space $L^0(\mathsf P)$ is continuous (see [1, chap.~5, \S\,3, Theorem 1]). Hence, the unit ball of $X$ is bounded in probability and for any $\beta>0$ there exists a sequence $M_n\uparrow\infty$ such that 
$$ \mathsf P(y\ge M_n)\le\frac{\beta}{2^{n+1}}\ \ \textrm{äëÿ âñåõ } y\in X_+,\ \|y\|\le 2.$$
Put
$$\varepsilon_n=\min\left\{\frac{1}{M_n},\frac{\beta}{2^n}\right\}.$$
Inequality (4) shows that there exists a sequence $z_n\in\conv |C_1|$ of the form
$$ z_n=x_n+y_n,\ \ x_n\in C_1,\ \ y_n\in X_+,\ \|y_n\|\le 2$$
where
$$ \mathsf P\left(A\bigcap\left\{x_n<\frac{1}{\varepsilon_n}\right\}\right)
<\frac{\varepsilon_n}{2}+\mathsf P\left(y_n\ge\frac{1}{\varepsilon_n}\right)\le\frac{\beta}{2^{n+1}}+\mathsf P(y_n\ge M_n)
\le \frac{\beta}{2^n}.$$

Let
$$ B_n=\{\omega\in A:\varepsilon_n x_n\ge 1\},\ \ \ B=\bigcap_{n=1}^\infty B_n.$$
Then $\mathsf P(A\backslash B_n)\le\beta/2^n$ and
$$ \mathsf P(A)=\mathsf P(B)+\mathsf P(A\backslash B)\le\mathsf P(B)+\sum_{n=1}^\infty \mathsf P(A\backslash B_n)\le\mathsf P(B)+\beta.$$
Choosing $\beta$ as any number in $(0,\mathsf P(A))$ we conclude that $\mathsf P(B)>0$.

Consider the sequence
$ v_n=I_B-\varepsilon_n x_n^-\in L^0(\mathsf P). $
The inequality
$$ |v_n|= I_B+ \varepsilon_n x_n^-\le \varepsilon_n x_n^++ \varepsilon_n x_n^-=\varepsilon_n |x_n|$$
implies that $v_n$ is an element of $X$. Moreover,
$$ \varepsilon_n x_n-v_n=\varepsilon_n x_n^+ - I_B\ge 0. $$
Thus, $v_n\in C-X_+=C$. From the closedness of the cone $C$ and the inequality
$$ \|v_n-I_B\|=\varepsilon_n\|x_n^-\|\le \varepsilon_n$$
it follows that $I_B\in C$. But this is impossible by virtue of (1). 

The obtained contradiction means that for any $A\in\mathscr F$ with $\mathsf P(A)>0$ there exists a number $c>0$ such that condition (2) is satisfied. As we mentioned above, this implies the assertion of Lemma 1. $\square$

It should be mentioned that in the course of the proof of Lemma 1 we took into account some ideas of [12], namely the notion of a hereditary unbounded set and the correspondent decomposition result ([12], Lemma 2.3).

In the next two lemmas we consider general Banach lattices [1, chap.10], [2, chap.9], which need not be ideal spaces. Denote by $U$, $U'$ the unit balls of Banach spaces $X$, $X'$ and put $U'_+=U'\cap X_+'$, where $X'_+=\{f\in X':\langle x,f\rangle\ge 0,\ x\in X_+\}$ is the non-negative cone of $X'$. Recall that the set
$H^\circ=\{f\in X':\langle x,f\rangle\le 1,\ x\in H\}$ is called a (one-sided) polar of $H\subset X$. The polar $W^\circ\subset X$ of a set $W\subset X'$ is defined similarly.

\begin{lemma}
For any element $x$ of a Banach lattice $X$ we have
$$ \|x^+\|=\sup\{\langle x,h\rangle:h\in U'_+\}.$$
\end{lemma}
\textsc{Proof.} Evidently, 
$$\sup\{\langle x,h\rangle:h\in U'_+\}\le \sup\{\langle x^+,h\rangle:h\in U'_+\}\le\|x^+\|.$$ 
Let us prove the converse inequality. If $\sup\{\langle x,h\rangle:h\in U'_+\}=0$, then $x\in (X'_+)^\circ=-X_+$.
Hence, $\|x^+\|=0.$ Assuming 
$$ \sup\{\langle x,h\rangle:h\in U'_+\}=\alpha>0$$
we shall check that $\|x^+\|\le\alpha$. 
Let
\begin{equation}
 \langle u-v,g\rangle\le 1\ \ \textrm{äëÿ âñåõ\ } u\in U, v\in X_+.
\end{equation} 
Putting $u=0$, we conclude that $g\in X'_+$. Together with inequality (5) this implies 
\begin{eqnarray*}
\|g\| &=& \sup\{|\langle w,g\rangle|:w\in U\}\le\sup\{\langle w^+,g\rangle+\langle w^-,g\rangle :w\in U\}\\
&=& \sup\{\langle |w|,g\rangle:w\in U\} 
      \le \sup\{\langle u,g\rangle:u\in U\}\le 1,
\end{eqnarray*}
that is, $g\in X'_+\cap U'=U'_+$. Conversely if $g\in U'_+$, then inequality (5) is satisfied. Thus, 
 $$ U'_+=(U-X_+)^\circ$$
and, by the bipolar theorem [2, Theorem 5.103], we have $(U'_+)^\circ=\overline{U-X_+}$, where the overline denotes the closure with respect to the norm of $X$. 

Furthermore, since $x/\alpha\in (U'_+)^\circ$, there exists a sequence $(y_n-z_n)_{n=1}^\infty$, $y_n\in U$, $z_n\in X_+$, converging to $x/\alpha$ in norm. This yields that
$$ \|x^+/\alpha\|=\lim_{n\to\infty}\|(y_n-z_n)^+\|\le\limsup_{n\to\infty}\|y_n^+\|\le 1.\ \ \square$$
\begin{lemma}
Let $X$ be a Banach lattice and let $C$ be a closed convex cone in $X$ with $C\cap X_+=\{0\}$. If there exists an element $f\in X'$ such that 
$$\sup_{x\in C_1}\langle x,f\rangle<\infty,  \ \ \ C _1=\{x\in C:\|x^-\|\le 1\}, $$
then there exists $g\in X'$ satisfying the conditions $f\le g$, $g\in C^\circ.$
\end{lemma}
\textsc{Proof.} Put $\lambda=\sup_{x\in C_1}\langle x,f\rangle$. If the assertion of lemma is false, then 
$$ (f+\lambda U'_+)\cap C^\circ=\varnothing.$$
By applying the separation theorem to the $\sigma (X',X)$-compact set $f+\lambda U'_+$ and 
$\sigma (X',X)$-closed set $C^\circ$, we conclude that there exists $x\in X$ such that
$$ \sup_{\eta\in C^\circ}\langle x,\eta\rangle<\inf\{\langle x,\zeta\rangle:\zeta\in f+\lambda U'_+\}.$$
Since $C$ is a cone, it follows that $\langle x,\eta\rangle\le 0$, $\eta\in C^\circ$ and $x\in C^{\circ\circ}=C$ by the bipolar theorem. Moreover,
\begin{eqnarray}
\langle x,f\rangle+\lambda\inf\{\langle x,h\rangle:h\in U'_+\}>0.
\end{eqnarray}
Note that
\begin{eqnarray}
\inf\{\langle x,h\rangle:h\in U'_+\}=-\sup\{\langle -x,h\rangle:h\in U'_+\}=-\|(-x)^+\|=-\|x^-\|
\end{eqnarray}
by Lemma 2.

By virtue of $0\neq x\in C$ we have $x\not\in X_+$ and $\|x^-\|>0$. From (6), (7) it follows that 
$$ \langle x/\|x^-\|,f\rangle>\lambda$$
in contradiction to the definition of $\lambda$ since $x/\|x^-\|\in C_1$. $\square$

For $X=L^\infty(\mathsf P)$ Lemma 3 was proved in [7, Lemma 2.5]. We also mention that this statement holds true for a non-closed cone $C$ as well. The correspondent result in the setting of locally convex-solid Riesz spaces was obtained in [13, Theorem 1].

\textsc{Proof of Theorem 1.} By virtue of Lemma 1 there exists a probability measure $\mathsf Q$, equivalent to $\mathsf P$, such that $C\subset L^1(\mathsf Q)$. It follows that $X\subset L^1(\mathsf Q)$ since $-X_+\subset C$ (compare with [14, Theorem 7]). Moreover, by Lemma 1 the strictly positive functional $x\mapsto\langle x,f\rangle=\mathsf E_\mathsf Q x$, $f\in L^\infty(\mathsf P)$, defined on $X$, is bounded from above on $C_1$. As is well-known, any positive functional on a Banach lattice is continuous [2, Theorem 9.6]. Therefore, $f\in X'$ and the desired element $g$ exists by Lemma 3.

\end{document}